\documentclass[11pt]{amsart}
\usepackage[margin=1.15in]{geometry}
\usepackage{amsmath,amssymb,amsthm,mathtools}
\usepackage{booktabs}
\usepackage{xcolor}
\usepackage[hypertexnames=false,colorlinks,linkcolor=blue!50!black,citecolor=green!35!black,urlcolor=blue!50!black]{hyperref}
\usepackage[capitalise,nameinlink]{cleveref}
\usepackage{thmtools}
\declaretheorem[style=plain,numberwithin=section]{theorem}
\declaretheorem[style=plain,sibling=theorem]{proposition,lemma,corollary}
\declaretheorem[style=plain,name=Question,numbered=no]{question*}
\declaretheorem[style=plain,name={Question~1$'$},numbered=no]{questionprime}
\declaretheorem[style=plain,name=Theorem,refname={Theorem,Theorems},Refname={Theorem,Theorems}]{introthm}

\declaretheorem[style=definition,sibling=theorem]{definition,example}
\declaretheorem[style=remark,sibling=theorem]{remark}
\numberwithin{equation}{section}

\newcommand{\hplus}{\boxplus}
\newcommand{\hminus}{\boxminus}
\newcommand{\val}{\operatorname{val}}
\newcommand{\m}{\mathfrak m}
\newcommand{\Hr}{\mathcal H_\rho}
\newcommand{\Mas}{\mathcal M}
\newcommand{\Lay}{\mathcal L}
\newcommand{\res}{\operatorname{res}}
\newcommand{\Z}{\mathbb Z}
\newcommand{\F}{\mathbb F}
\newcommand{\cls}[1]{[#1]}
\newcommand{\lt}[2]{#1\,t^{#2}}

\begin{document}

\title{The norm axiom and factor hyperfields}

\author{Alessandro Linzi}
\address{Independent researcher, Udine, Italy}
\email{alessandro.linzi.phd@icloud.com}


\begin{abstract}
Touchard and the author asked whether every valued hyperfield, in Krasner's sense with the norm axiom, is of the form $\mathcal H_\rho(K)$: a valued field $K$ read up to relative error $\rho$. The question turns on the norm axiom. If a valued hyperfield is a quotient $K^\times/T\cup\{0\}$ of a field, the norm axiom forces $T=1+\mathfrak m_\rho$ for a valuation on $K$. This completes the proof of a published proposition, and it makes the question equivalent to asking whether every Krasner valued hyperfield is such a quotient. For value group $\mathbb Z$ we prove directly that it is. Without the norm axiom the answer is negative. Layering a hyperfield over a non-trivial ordered abelian group gives a factor hyperfield exactly when the hyperfield itself is one, so Massouros's hyperfields yield hyperfields with a non-trivial valuation, satisfying all the other axioms, that are not factor hyperfields. The general question remains open.
\end{abstract}

\maketitle

\section{Introduction}\label{sec:intro}

\subsection{The question}
Let $K$ be a field with a valuation $v$ and value group $\Gamma$. Fix an initial segment $\rho$ of $\Gamma_{\ge 0}$ and let $\m_\rho$ be the set of elements of $K$ whose value exceeds every element of $\rho$. The elements of $1+\m_\rho$ form a subgroup of $K^\times$. The quotient
\[
  \Hr(K) \coloneqq K^\times/(1+\m_\rho)\cup\{0\}
\]
remembers each element of $K$ up to a relative error of size $\rho$. Multiplication passes to the quotient, but addition does not: the sum of two classes is the \emph{set} of classes of all possible sums of representatives. This makes $\Hr(K)$ a \emph{hyperfield}, and the valuation of $K$ descends to it. For $\rho=\{0\}$ one obtains the leading-term structure $\mathrm{RV}(K)$ used in the model theory of valued fields; for discrete $\Gamma$ and $\rho=\{0,\dots,n-1\}$ one obtains the hyperfields studied by Krasner~\cite{Krasner57}, Tolliver~\cite{Tolliver16} and Lee~\cite{Lee20}.

Quotients of a field, or more generally of a ring, by a multiplicative subgroup are a classical source of hyperstructures, going back to Krasner~\cite{Krasner83}. The Krasner hyperfield and the hyperfield of signs, which Baker and Bowler use as coefficients for matroids and oriented matroids, arise in this way \cite[Remark~2.7 and Examples~2.10, 2.13]{BakerBowler19}, and so does the hyperring structure that Connes and Consani put on the ad\`ele class space of a global field \cite{ConnesConsani11}.

Krasner isolated the properties of $\Hr(K)$ that make it a \emph{valued hyperfield}. The decisive one is the \emph{norm axiom}: each sum $x\hplus y$ is a full ultrametric ball of a prescribed radius. Touchard and the author~\cite{LT25} studied when systems of valued hyperfields come from a single valued field, and at the end of their paper they asked the converse question for a single hyperfield. We quote it verbatim.

\begin{question*}[{\cite[Question~1, p.~23]{LT25}}]
Are all valued hyperfields with value group $\Gamma$ of the form $\mathcal H_\rho(K)$ for some valued field $K$ and some initial segment $\rho$ of $\Gamma_{\ge 0}$?
\end{question*}

Here ``valued hyperfield'' means a hyperfield with a valuation of norm $\rho$ in the sense of \cite[Definition~4.6]{Linzi25}, recalled in \cref{def:valuation}; by \cite[Proposition~4.8(i)]{Linzi25} this is equivalent to \cite[Definition~1.8]{LT25}, as explained after that definition. The same circle of ideas appears in two places in~\cite{Linzi25}. Proposition~4.11 there states that a factor hyperfield carrying a Krasner valuation is of the form $\Hr(K)$, and Remark~4.12 records that it is not known whether every Krasner valued hyperfield is a factor hyperfield. Section~7, item~(2), asks the analogous question for the weaker notion of valuation of \cite[Definition~3.1]{Linzi25}, which consists of the axioms (V0)--(V2) below (labelled (V1)--(V3) there): \emph{is it true that any hyperfield admitting a non-trivial valuation is a factor hyperfield?}

\subsection{Results}
The three results below show that the question is about what the norm axiom (V3) forces. With (V3), a presentation as a factor hyperfield is automatically a presentation as $\Hr(K)$ (\cref{thm:A}), and in the discrete case the axioms rebuild the field without any factor presentation being given (\cref{thm:Z}). Without (V3), a non-trivial valuation does not force a factor presentation (\cref{thm:transfer}).

A \emph{factor hyperfield} is a quotient $K_T=K^\times/T\cup\{0\}$ of a field $K$ by an arbitrary subgroup $T\le K^\times$, with the induced hyperaddition (\cref{sec:prelim}). Every $\Hr(K)$ is a factor hyperfield, with $T=1+\m_\rho$. The first result is that the converse holds as soon as a Krasner valuation is present.

\begin{introthm}[{completing \cite[Proposition~4.11]{Linzi25}}]\label{thm:A}
Let $(H,\val)$ be a valued hyperfield of norm $\rho$ and value group $\Gamma$. Suppose that $H$ is isomorphic, as a hyperfield, to a factor hyperfield $K_T$. Then $K$ carries a valuation $w$ with value group $\Gamma$ such that $T=1+\m_\rho(w)$, and the isomorphism $K_T\to H$ is an isomorphism of valued hyperfields $\Hr(K,w)\to(H,\val)$.
\end{introthm}

This is the statement of \cite[Proposition~4.11]{Linzi25}, whose proof in \cite{Linzi25} shows only that $T$ consists of $1$-units, that is, $T\subseteq1+\m_{\{0\}}(w)$. The missing equality $T=1+\m_\rho(w)$ is Krasner's observation \cite[p.~148]{Krasner57} in the form of \cite[Lemma~1.25]{LT25}; \cref{sec:A} spells out the combination. As a consequence, the question of Touchard and the author is equivalent to the question of \cite[Remark~4.12]{Linzi25}: it has a positive answer if and only if every Krasner valued hyperfield is a factor hyperfield (\cref{cor:equivalent}).

The second result settles the question for the value group $\Z$.

\begin{introthm}\label{thm:Z}
Let $(H,\val)$ be a valued hyperfield of norm $\rho$ with value group $\Z$. Then $H\cong\Hr(K)$ as valued hyperfields for a discretely valued field $K$.
\end{introthm}

This can also be derived from Tolliver's comparison of Krasner's discretely valued hyperfields with Deligne's triples~\cite{Tolliver16}, together with \cref{thm:A}; see \cref{rem:tolliver}. Our proof follows the same route, through the ring of \cite[Lemmas~4.2 and~4.3]{Tolliver16} and Deligne's observation that truncated discrete valuation rings are quotients of discrete valuation rings, but it is written directly in the axioms of~\cite{LT25} and is self-contained apart from a theorem of Hungerford~\cite{Hungerford68}. Lee \cite[Theorems~4.8 and~4.10]{Lee20} proves much more in a narrower setting: for discretely valued hyperfields of mixed characteristic $(0,p)$ with perfect residue field and length $n$ large compared with the ramification index, he obtains equivalences of categories with truncated discrete valuation rings and with complete discrete valuation rings. \Cref{thm:Z} needs none of these hypotheses, but it is only an existence statement.

The third result answers the question of \cite[\S7(2)]{Linzi25} in the negative. It uses the layering $\Mas\rtimes\Gamma$ of a hyperfield $\Mas$ along an ordered abelian group, due to Bowler and Su~\cite[Definition~4.1]{BowlerSu21}: one copy of $\Mas$ is placed in each degree $\gamma\in\Gamma$, as residue field and value group are stacked in a power series field (\cref{sec:layer}).

\begin{introthm}\label{thm:transfer}
Let $\Mas$ be a hyperfield and $\Gamma$ a non-trivial ordered abelian group. The layering $\Mas\rtimes\Gamma$ carries a surjective map onto $\Gamma\cup\{\infty\}$ satisfying the axioms \textup{(V0)--(V2)} of a valuation, and $\Mas\rtimes\Gamma$ is a factor hyperfield if and only if $\Mas$ is. In particular, if $\Mas=\Mas(G)$ is a Massouros hyperfield with $G$ of exponent~$2$ and $|G|\ge4$, then $\Mas\rtimes\Gamma$ is not a factor hyperfield; it admits no valuation of any norm in the sense of \cref{def:valuation}.
\end{introthm}

Massouros~\cite{Massouros85} constructed the hyperfields $\Mas(G)$ and showed that they are not factor hyperfields; as observed in \cite[\S7(2)]{Linzi25}, they admit only the trivial valuation. One direction of the equivalence is the residue computation of \cite[Lemma~3.6 and Proposition~3.14]{Linzi25}. The other uses the field of fractions of the group algebra $k[t^\Gamma]$, in which a factor presentation $k_S$ of $\Mas$ propagates to one of $\Mas\rtimes\Gamma$. For the Krasner hyperfield, the hyperfield of signs and fields, that is, for $S=k^\times$, $S$ the positive cone of an ordered field, and $S=\{1\}$, this construction is due to Bowler and Su \cite[proof of Theorem~7.5, Cases~1--3]{BowlerSu21}, who use Hahn-type fields; the Krasner case is also \cite[Proposition~2.17(iv)]{Linzi25}. \Cref{prop:lift} extends it to an arbitrary subgroup $S$. The norm axiom fails because $1\hminus1$ contains elements of value $0$, the mechanism behind \cite[Proposition~4.14]{Linzi25}. \Cref{tab:examples} shows four instances.

\begin{table}[htbp]
\centering
\caption{Layerings $\Mas\rtimes\Gamma$ over a non-trivial ordered abelian group $\Gamma$. ``Valued'' means that $\Mas\rtimes\Gamma$ admits a valuation of some norm in the sense of \cref{def:valuation}. The third column follows from \cref{thm:transfer}; for rows 1--3 it was known \cite[proof of Theorem~7.5, Cases~1--3]{BowlerSu21}, \cite[Proposition~2.17(iv)]{Linzi25}. The fourth column comes from the paragraph after \cref{prop:lift} (row 1), \cref{rem:table} (rows 2 and 3) and \cref{lem:not-krasner} (row 4).}\label{tab:examples}
\begin{tabular}{@{}llcc@{}}
\toprule
$\Mas$ & $\Mas\rtimes\Gamma$ & factor? & valued? \\
\midrule
a field $k$ & $\mathcal H_{\{0\}}(k(t^\Gamma))$ & yes & yes \\
Krasner hyperfield & tropical hyperfield of $\Gamma$ & yes & no \\
hyperfield of signs & tropical real hyperfield of $\Gamma$ & yes & no \\
$\Mas(G)$, $G^2=1$, $|G|\ge4$ & \cref{thm:transfer} & no & no \\
\bottomrule
\end{tabular}
\end{table}

For $\Gamma=\mathbb R$ the second and third rows are Viro's tropical hyperfield and tropical real hyperfield \cite[\S\S5.3 and~7.2]{Viro10}; see \cref{rem:table}.

\subsection{What is new}
The mathematics in this note is elementary, and several steps are in the literature. What we add is a complete proof of \cite[Proposition~4.11]{Linzi25} and the resulting equivalence of the two published questions (\cref{thm:A,cor:equivalent}); a proof in the language of \cite{LT25} that the question of Touchard and the author has a positive answer for the value group $\Z$ (\cref{thm:Z}); and the transfer principle of \cref{thm:transfer}, which turns every hyperfield that is not a factor hyperfield into a negative answer to \cite[\S7(2)]{Linzi25}. For the Krasner hyperfield, the hyperfield of signs and fields, the fact that the layering is a factor hyperfield is due to Bowler and Su \cite[proof of Theorem~7.5, Cases~1--3]{BowlerSu21} and, in the Krasner case, to \cite[Proposition~2.17(iv)]{Linzi25}; \cref{prop:lift} extends their construction to every subgroup $S$. New are the equivalence for every hyperfield $\Mas$, the direction from $\Mas\rtimes\Gamma$ to $\Mas$ written out in full, and the application to Massouros hyperfields.

\subsection{What remains open}
The main theorem of Touchard and the author \cite[Theorem~4.1]{LT25} gives a positive answer when the norm $\rho$ is the non-negative part $\Delta_{\ge0}$ of a convex subgroup $\Delta$ of $\Gamma$, for instance $\rho=\{0\}$ or $\rho=\Gamma_{\ge0}$ (\cref{cor:convex}). With \cref{thm:Z}, this leaves value groups other than $\Z$ with norms not of the form $\Delta_{\ge0}$. The following is the most concrete such case that we cannot settle.

\begin{questionprime}
Let $\Gamma=\Z\times\Z$ with the lexicographic order and $\rho=\{(0,0),(0,1)\}$. Is every valued hyperfield of norm $\rho$ with value group $\Gamma$ of the form $\Hr(K)$?
\end{questionprime}

The convex subgroups of $\Gamma$ are $\{0\}$, $\{0\}\times\Z$ and $\Gamma$, with non-negative parts $\{(0,0)\}$, $\{(0,m):m\ge0\}$ and $\Gamma_{\ge0}$, so $\rho$ is not of the form $\Delta_{\ge0}$ and \cref{cor:convex} does not apply; \cref{thm:Z} does not apply because $\Gamma\not\cong\Z$. See \cref{rem:heuristic} for a heuristic.

\subsection{Organisation}
\Cref{sec:prelim} recalls hyperfields, factor hyperfields and Krasner valuations. \Cref{sec:A} proves \cref{thm:A} and records what \cite[Theorem~4.1]{LT25} gives. \Cref{sec:C} treats the value group $\Z$. \Cref{sec:massouros,sec:layer,sec:B} contain Massouros hyperfields, layerings and the proof of \cref{thm:transfer}. \Cref{sec:checks} describes the computer checks that accompany the paper.

\subsection*{Running example}
Throughout, $\Mas_4$ denotes Massouros's hyperfield on the Klein four-group $V=\{1,a,b,ab\}$ together with $0$, and $\Lay_4$ denotes its layering over $\Z$. Both are small enough to compute with by hand.

\section{Preliminaries}\label{sec:prelim}

We follow \cite{LT25} closely.

\subsection{Hyperfields}
A \emph{hyperoperation} on a set $H$ assigns to each pair $x,y\in H$ a non-empty subset $x\hplus y\subseteq H$. For subsets $A,B$ we put $A\hplus B=\bigcup_{a\in A,\,b\in B}a\hplus b$, and $x\hplus A=\{x\}\hplus A$, $A\hplus x=A\hplus\{x\}$.

\begin{definition}[{\cite[Definitions~1.1 and~1.4]{LT25}}]\label{def:hyperfield}
A \emph{hyperfield} is a tuple $(H,\hplus,\cdot,0,1)$ such that:
\begin{enumerate}
\item[(CH1)] $(x\hplus y)\hplus z=x\hplus(y\hplus z)$ for all $x,y,z$;
\item[(CH2)] $x\hplus y=y\hplus x$ for all $x,y$;
\item[(CH3)] for each $x$ there is a unique $-x$ with $0\in x\hplus(-x)$;
\item[(CH4)] if $z\in x\hplus y$ then $x\in z\hminus y$;
\item[(HF2)] $H^\times=H\setminus\{0\}$ is an abelian group under $\cdot$ and $x\cdot0=0$ for all $x$;
\item[(HF3)] $z\cdot(x\hplus y)=zx\hplus zy$ for all $x,y,z$.
\end{enumerate}
\end{definition}

These axioms imply $x\hplus0=\{x\}$ \cite[Remark~1.2]{LT25}. A \emph{homomorphism} $\varphi\colon H\to H'$ of hyperfields satisfies $\varphi(0)=0$, $\varphi(xy)=\varphi(x)\varphi(y)$ and $\varphi(x\hplus y)\subseteq\varphi(x)\hplus\varphi(y)$; an \emph{isomorphism} is a bijection $\varphi$ with $\varphi(x\hplus y)=\varphi(x)\hplus\varphi(y)$ and $\varphi(xy)=\varphi(x)\varphi(y)$. An isomorphism maps $1$ to $1$ and $-x$ to $-\varphi(x)$, by uniqueness in (HF2) and (CH3).

\subsection{Factor hyperfields}
Let $K$ be a field and $T$ a subgroup of $K^\times$. Write $\cls x_T\coloneqq xT$ for $x\in K^\times$ and $\cls0_T\coloneqq\{0\}$. On $K_T\coloneqq K^\times/T\cup\{\cls0_T\}$ put, for $x,y\in K$,
\begin{equation}\label{eq:factor}
  \cls x_T\hplus\cls y_T\coloneqq\{\cls{x+ys}_T : s\in T\},\qquad \cls x_T\cdot\cls y_T\coloneqq\cls{xy}_T.
\end{equation}
Krasner showed that $K_T$ is a hyperfield~\cite{Krasner83}; see also \cite[Proposition~1.15]{LT25}. Hyperfields isomorphic to some $K_T$ are called \emph{factor} (or \emph{quotient}) hyperfields. Two facts are used repeatedly: $-\cls x_T=\cls{-x}_T$, and $\cls{x+y}_T\in\cls x_T\hplus\cls y_T$ for all $x,y\in K$ (take $s=1$).

\subsection{Valuations}
\begin{definition}[{\cite[Definition~4.6]{Linzi25}; equivalent to \cite[Definition~1.8]{LT25}}]\label{def:valuation}
Let $\Gamma$ be an ordered abelian group and $\rho$ a non-empty initial segment of $\Gamma_{\ge0}$, that is, a non-empty subset of $\Gamma_{\ge0}$ closed downwards in $\Gamma_{\ge0}$; in particular $0\in\rho$. As usual, $\infty$ is a new symbol with $\gamma<\infty$ and $\gamma+\infty=\infty+\gamma=\infty$ for all $\gamma\in\Gamma$. A \emph{valuation of norm $\rho$} on a hyperfield $H$ is a surjective map $\val\colon H\to\Gamma\cup\{\infty\}$ such that, for all $x,y\in H$:
\begin{enumerate}
\item[(V0)] $\val(x)=\infty$ if and only if $x=0$;
\item[(V1)] $\val(xy)=\val(x)+\val(y)$;
\item[(V2)] $\val(z)\ge\min\{\val(x),\val(y)\}$ for all $z\in x\hplus y$;
\item[(V3)] \textup{(norm axiom)} if $z\in x\hplus y$, then for all $z'\in H$,
\[
 z'\in x\hplus y \iff \val(u)>\rho+\min\{\val(x),\val(y)\}\ \text{for all } u\in z\hminus z'.
\]
\end{enumerate}
A \emph{valued hyperfield} is a hyperfield with a valuation of some norm; in \cite[Definition~4.6]{Linzi25} these are called Krasner valued hyperfields. An \emph{isomorphism of valued hyperfields} $(H,\val)\to(H',\val')$ is an isomorphism of hyperfields $\varphi$ with $\val'\circ\varphi=\val$.
\end{definition}

In (V3), $\rho+\gamma$ is the smallest initial segment of $\Gamma$ containing $\delta+\gamma$ for all $\delta\in\rho$, and $\val(u)>\rho+\gamma$ means $\val(u)-\gamma>\rho$ (true for $u=0$). In a valued field or a valued hyperfield we write $\m_\rho\coloneqq\{x:\val(x)>\rho\}$; this set contains $0$. When the valuation must be named we write $\m_\rho(w)$. The \emph{trivial} valuation is the one with $\Gamma=\{0\}$. A map satisfying only (V0)--(V2) is called a valuation in \cite[Definition~4.1]{KLS22} and \cite[Definition~3.1]{Linzi25}; see \cite[Remark~1.9]{LT25}. \Cref{def:valuation} is \cite[Definition~4.6]{Linzi25}: such a valuation together with the axiom (KVH) there, which is (V3). The only difference of presentation is that the norm $\rho_v$ there is an initial segment of the value group containing $0$, and here we take $\rho=\rho_v\cap\Gamma_{\ge0}$. By \cite[Proposition~4.8(i)]{Linzi25} (replace $y$ by $-y$), this is equivalent to \cite[Definition~1.8]{LT25}, which includes in addition the single-value property that if $0\notin x\hplus y$ then $\{\val(z):z\in x\hplus y\}$ is a singleton. Thus the Question and the results of~\cite{LT25} used below apply.

If $(K,v)$ is a valued field and $T=1+\m_\rho$, then the elements of $T$ have value $0$, so $v$ descends to a map $v_\rho$ on $K_T$; it is a valuation of norm $\rho$ \cite[Section~1.2]{LT25}. We write $\Hr(K)$ or $\Hr(K,v)$ for this valued hyperfield.

\begin{example}\label{ex:Hrho}
Take $K=\mathbb Q$ with the $3$-adic valuation and $\rho=\{0\}$. Then $1+\m_\rho=1+3\Z_{(3)}$, so a class in $\Hr(K)$ records the exact power of $3$ and the leading $3$-adic digit, which is $1$ or $2$. In it, $\cls1\hplus\cls1=\{\cls2\}$, while $\cls1\hplus\cls{2}$ consists of $0$ and all classes of positive value, the open ball of radius $0$ around $0$.
\end{example}

Two consequences of the norm axiom will be used repeatedly. Both are contained in the proof of \cite[Lemma~1.25]{LT25}, and the first is also the key step of \cite[Proposition~4.14]{Linzi25}. The first says that $1\hminus1$ is exactly the open ball of radius $\rho$ around $0$.

\begin{lemma}\label{lem:one-minus-one}
Let $(H,\val)$ be a valued hyperfield of norm $\rho$. Then $1\hminus 1=\m_\rho$.
\end{lemma}

\begin{proof}
We have $0\in 1\hminus 1$, and $\min\{\val(1),\val(-1)\}=0$ by \cite[Lemma~1.10]{LT25}. Apply (V3) with $x=1$, $y=-1$ and $z=0$: an element $z'$ lies in $1\hminus1$ if and only if every $u\in 0\hminus z'$ satisfies $\val(u)>\rho$. Since $0\hminus z'=\{-z'\}$ and $\val(-z')=\val(z')$, this says $\val(z')>\rho$.
\end{proof}

The second says that adding something $\rho$-small to $1$ changes nothing.

\begin{lemma}\label{lem:absorb}
Let $(H,\val)$ be a valued hyperfield of norm $\rho$. If $\val(m)>\rho$, then $1\hplus m=\{1\}$.
\end{lemma}

\begin{proof}
By \cref{lem:one-minus-one}, $m\in 1\hplus(-1)$, so $1\in m\hplus 1$ by (CH4) and (CH2). By \cite[Lemma~1.12]{LT25}, $1\in 1\hplus m$ forces $1\hplus m=\{1\}$.
\end{proof}

\begin{remark}[The norm is determined by the valuation]\label{rem:norm-intrinsic}
Let $(H,\val)$ be a valued hyperfield of norm $\rho$. Then
\[
  \rho=\Gamma_{\ge0}\setminus\val(1\hminus1).
\]
Indeed, by \cref{lem:one-minus-one} and the surjectivity of $\val$, $\val(1\hminus1)=\{\gamma\in\Gamma:\gamma>\rho\}\cup\{\infty\}$. An element $\gamma\ge0$ fails to be $>\rho$ exactly when $\gamma\le\delta$ for some $\delta\in\rho$, that is, exactly when $\gamma\in\rho$, because $\rho$ is an initial segment of $\Gamma_{\ge0}$. So a map $\val$ is a valuation of at most one norm, and the norm is not extra data. What can vary is the valuation itself: a hyperfield can carry inequivalent Krasner valuations \cite[Example~4.24]{Linzi25}.
\end{remark}

\section{Valued factor hyperfields: proof of Theorem~\ref{thm:A}}\label{sec:A}

\Cref{thm:A} combines two published facts. By \cite[Lemma~3.6]{Linzi25}, a valuation on a factor hyperfield $K_T$ pulls back to a valuation on $K$. The resulting map from $K$ onto $H$ is then a surjective isometric homomorphism of valued hyperfields in the sense of \cite[Definition~1.18]{LT25}, and Krasner's observation \cite[p.~148]{Krasner57}, as formulated in \cite[Lemma~1.25]{LT25}, says that the kernel of such a map is the group of $1$-units of the right level. The proof of \cite[Proposition~4.11]{Linzi25} contains the first step and the inclusion $T\subseteq 1+\m_{\{0\}}$, but not the equality $T=1+\m_\rho$. For the reader's convenience we give the short argument in full. The norm axiom pins down $T$ from both sides: \cref{lem:one-minus-one} gives $T\subseteq 1+\m_\rho$ and \cref{lem:absorb} gives $1+\m_\rho\subseteq T$.

\begin{proof}[Proof of \cref{thm:A}]
Let $\varphi\colon K_T\to H$ be an isomorphism of hyperfields and define $w\colon K\to\Gamma\cup\{\infty\}$ by $w(x)\coloneqq\val(\varphi(\cls x_T))$.

\emph{Step 1: $w$ is a valuation on the field $K$ with value group $\Gamma$} \cite[Lemma~3.6]{Linzi25}. By (V0), $w(x)=\infty$ exactly when $x=0$. Multiplicativity follows from (V1). For the ultrametric inequality, $\cls{x+y}_T\in\cls x_T\hplus\cls y_T$, hence $\varphi(\cls{x+y}_T)\in\varphi(\cls x_T)\hplus\varphi(\cls y_T)$, and (V2) gives $w(x+y)\ge\min\{w(x),w(y)\}$. Since $\varphi$ is a bijection and $\val$ is surjective, $w(K^\times)=\Gamma$.

\emph{Step 2: $T\subseteq 1+\m_\rho(w)$.} This is clear for $s=1$. Let $s\in T$ with $s\ne1$. By \eqref{eq:factor}, $\cls{1-s}_T\in\cls1_T\hplus\cls{-1}_T$, so $\varphi(\cls{1-s}_T)\in 1\hminus1$. By \cref{lem:one-minus-one}, $w(1-s)>\rho$.

\emph{Step 3: $1+\m_\rho(w)\subseteq T$.} This is clear for $m=0$. Let $m\in K^\times$ with $w(m)>\rho$. By \cref{lem:absorb}, $1\hplus\varphi(\cls m_T)=\{1\}$, so $\cls1_T\hplus\cls m_T=\{\cls1_T\}$. In particular $\cls{1+m}_T=\cls1_T$, that is, $1+m\in T$.

\emph{Step 4: conclusion.} By Steps 2 and 3, $T=1+\m_\rho(w)$, so $K_T=\Hr(K,w)$ as hyperfields. The valuation of $\Hr(K,w)$ sends $\cls x_T$ to $w(x)=\val(\varphi(\cls x_T))$, so $\varphi$ preserves valuations.
\end{proof}

\begin{remark}[Where the norm axiom is used]\label{rem:audit}
Step~1 uses only (V0)--(V2). Step~2 uses (V3) through the inclusion $1\hminus1\subseteq\m_\rho$ of \cref{lem:one-minus-one}, and Step~3 uses it through \cref{lem:absorb}, at the given norm $\rho$. Two factor hyperfields with maps satisfying (V0)--(V2) show that Steps~2 and~3 can each fail without (V3) at the given norm.
\begin{enumerate}
\item $\F_5/\{\pm1\}$ with the trivial valuation. The only valuation on $\F_5$ is trivial, so $1+\m_\rho(w)=\{1\}$ for every $\rho$, and Step~3 holds trivially. Step~2 fails: $T=\{\pm1\}$, and $1\hminus1=\{0,\cls2\}$ contains an element of value $0$.
\item $\mathcal H_{\{0,1\}}(\mathbb Q_p)$ with the $p$-adic valuation, regarded as a candidate valuation of norm $\{0\}$. It satisfies (V0)--(V2), and (V3) for the norm $\{0,1\}$ but not for $\{0\}$, in accordance with \cref{rem:norm-intrinsic}. Step~2 holds, since $T=1+p^2\Z_p\subseteq1+p\Z_p=1+\m_{\{0\}}$, but Step~3 fails: $1+p\notin T$, and $1\hplus\cls p\ne\{1\}$.
\end{enumerate}
\end{remark}

\begin{corollary}\label{cor:equivalent}
The following are equivalent.
\begin{enumerate}
\item The question of Touchard and the author has a positive answer: every valued hyperfield $(H,\val)$ of norm $\rho$ and value group $\Gamma$ is isomorphic, as a valued hyperfield, to $\Hr(K)$ for some valued field $K$ with value group $\Gamma$.
\item Every hyperfield that admits a valuation in the sense of \cref{def:valuation} is a factor hyperfield.
\end{enumerate}
\end{corollary}

Statement (2) is the open question of \cite[Remark~4.12]{Linzi25}. The question of Touchard and the author as quoted, which allows any $K$ and any $\rho$, lies between (1) and (2), so it is equivalent to both.

\begin{proof}
Each $\Hr(K)$ is a factor hyperfield, so (1) implies (2). Conversely, (2) and \cref{thm:A} give (1), with the same norm $\rho$.
\end{proof}

\begin{remark}
\Cref{cor:equivalent} also settles a point of interpretation: it does not matter whether the question asks for an isomorphism of hyperfields or of valued hyperfields. Indeed, if $H\cong\Hr(K)$ merely as hyperfields, then \cref{thm:A} provides a valuation $w$ on $K$, possibly different from the given one, with $H\cong\Hr(K,w)$ as valued hyperfields. Some care is needed here because a hyperfield can carry inequivalent Krasner valuations \cite[Example~4.24]{Linzi25}; \cref{thm:A} chooses $w$ to match the given valuation of $H$.
\end{remark}

\begin{remark}[Rigidity]\label{rem:rigid}
The norm axiom also determines a valued hyperaddition from below. Let $H$ carry two hyperadditions $+$ and $\hplus$, one multiplication and one map $\val$, such that $(H,+,\cdot,\val)$ and $(H,\hplus,\cdot,\val)$ are both valued hyperfields of norm $\rho$. If $x\hplus y\subseteq x+y$ for all $x,y$, then the two hyperadditions coincide. Indeed, $0\in x\hplus y$ implies $0\in x+y$, so both have the same negation and $z\hminus z'\subseteq z-z'$. Given $x,y$, pick $z\in x\hplus y\subseteq x+y$ and put $r=\rho+\min\{\val(x),\val(y)\}$. By (V3) in both structures, $x+y=\{z':\val(u)>r\ \forall u\in z-z'\}$ and $x\hplus y=\{z':\val(u)>r\ \forall u\in z\hminus z'\}$; the second condition quantifies over a smaller set, so $x+y\subseteq x\hplus y$. In particular the sums of $\Hr(K)$ cannot be replaced by strictly smaller ones without losing the norm axiom.
\end{remark}

\subsection*{Norms that are positive parts of convex subgroups}
The main theorem of Touchard and the author \cite[Theorem~4.1]{LT25} constructs a valued field from an isometric system of valued hyperfields whose norms exhaust the non-negative part of a convex subgroup. A single valued hyperfield is such a system with one member, so we obtain:

\begin{corollary}[{of \cite[Theorem~4.1]{LT25}}]\label{cor:convex}
Let $(H,\val)$ be a valued hyperfield of norm $\rho$ and value group $\Gamma$. If $\rho=\Delta_{\ge0}$ for a convex subgroup $\Delta$ of $\Gamma$, then $H\cong\Hr(K)$ as valued hyperfields for some valued field $K$ with value group $\Gamma$.
\end{corollary}

The statement of \cite[Theorem~4.1]{LT25} only asserts an isomorphism $\Hr(K)\simeq H$; that it can be taken to preserve valuations is visible in its proof, and it also follows from \cref{thm:A}.

This covers $\rho=\{0\}$ (the case $\Delta=\{0\}$, where $H$ is an $\mathrm{RV}$-structure) and $\rho=\Gamma_{\ge0}$ (where $H$ is a valued field). The question is therefore open only for norms that are not of this form, and by \cref{thm:Z} only for value groups other than $\Z$.

\begin{remark}\label{rem:leading-layer}
Let $(H,\val)$ be any valued hyperfield of norm $\rho$. By \cite[Proposition~1.17]{LT25}, the set $T_0\coloneqq 1\hplus\m_{\{0\}}$ is a subgroup of $H^\times$, and the quotient $H_{\{0\}}\coloneqq H_{T_0}$ is a valued hyperfield of norm $\{0\}$. By \cref{cor:convex} it is of the form $\mathcal H_{\{0\}}(K)$. So the leading-term layer of every valued hyperfield comes from a valued field, and a counterexample to the question would have to fail at a finer level of precision.
\end{remark}

\section{The discrete case: proof of Theorem~\ref{thm:Z}}\label{sec:C}

The idea is the one behind Deligne's triples~\cite{Deligne84}, as related to Krasner's hyperfields by Tolliver~\cite{Tolliver16}: inside a valued hyperfield $H$ of norm $\{0,\dots,n-1\}$, the elements of non-negative value, read modulo the open ball of radius $n-1$, form a ring $R$; this ring is a quotient of a discrete valuation ring $\mathcal O$, and $H$ is recovered as $\Hr(\operatorname{Frac}\mathcal O)$. In $\Hr(K)$ itself, $R$ is just $\mathcal O_K/\m_K^n$. \Cref{lem:sim,lem:R} correspond to \cite[Lemmas~4.2 and~4.3]{Tolliver16} (for the ring $M_0$ there), restated and reproved in the axioms of \cref{def:valuation}; part~(2) of \cref{lem:sim} is Krasner's \cite[p.~145]{Krasner57}, used in the proof of \cite[Proposition~6.1]{Tolliver16}.

Throughout this section $(H,\val)$ is a valued hyperfield with value group $\Z$ and norm $\rho$. Every initial segment of $\Z_{\ge0}$ is $\Z_{\ge0}$ or $\{0,\dots,n-1\}$ for some $n\ge1$. The norms $\Z_{\ge0}$ and $\{0\}$ are covered by \cref{cor:convex}, with $\Delta=\Z$ and $\Delta=\{0\}$, so we assume $\rho=\{0,\dots,n-1\}$ with $n\ge2$. Then $\val(u)>\rho+\gamma$ means $\val(u)\ge n+\gamma$. Put $\mathcal O_H\coloneqq\{x\in H:\val(x)\ge0\}$ and fix $\varpi\in H$ with $\val(\varpi)=1$ (in $\Hr(K)$, the class of a uniformiser).

\begin{lemma}\label{lem:sim}
For $x,y\in\mathcal O_H$ write $x\sim y$ if $\val(u)\ge n$ for all $u\in x\hminus y$.
\begin{enumerate}
\item $\sim$ is an equivalence relation on $\mathcal O_H$, and $x\sim0$ if and only if $\val(x)\ge n$.
\item If $\val(x)=0$, the class of $x$ is $\{x\}$.
\item If $x,y\in\mathcal O_H$ and $\min\{\val(x),\val(y)\}=0$, then $x\hplus y$ is a single $\sim$-class.
\item If $x\sim x'$, $y\sim y'$, $z\in x\hplus y$ and $z'\in x'\hplus y'$, then $z\sim z'$ and $xy\sim x'y'$.
\end{enumerate}
\end{lemma}

\begin{proof}
(1) Reflexivity: $0\hminus0=\{0\}$, and for $x\ne0$ we have $x\hminus x=x\,\m_\rho$ by \cref{lem:one-minus-one} and (HF3), whose elements have value $>\rho+\val(x)$, so value $\ge n$. Symmetry: $y\hminus x=-(x\hminus y)$. Transitivity: since $0\in y\hminus y$, associativity and commutativity give $x\hminus z\subseteq x\hplus(y\hminus y)\hminus z=(x\hminus y)\hplus(y\hminus z)$, and by (V2) a sum of two elements of value $\ge n$ has only elements of value $\ge n$. Finally $x\hminus0=\{x\}$.

(2) Apply (V3) with the pair $(x,0)$ and $z=x\in x\hplus0$: an element $z'$ lies in $x\hplus0=\{x\}$ if and only if every $u\in x\hminus z'$ has value $>\rho$, that is, if and only if $z'\sim x$.

(3) By (V2), $x\hplus y\subseteq\mathcal O_H$. Pick $z\in x\hplus y$. By (V3), an element $z'\in H$ lies in $x\hplus y$ if and only if every $u\in z\hminus z'$ has value $>\rho+0$. Hence $x\hplus y$ is the set of $z'\in\mathcal O_H$ with $z'\sim z$.

(4) As in (1), $z\hminus z'\subseteq(x\hplus y)\hminus(x'\hplus y')=(x\hminus x')\hplus(y\hminus y')$, all of whose elements have value $\ge n$. For the product, $xy\hminus x'y'\subseteq (x\hminus x')y\hplus x'(y\hminus y')$ by the same trick with $0\in x'y\hminus x'y$, and both summands consist of elements of value $\ge n$.
\end{proof}

Let $R\coloneqq\mathcal O_H/{\sim}$ and write $\bar x$ for the class of $x$. By \cref{lem:sim}(4), $\bar x+\bar y\coloneqq\bar z$ for any $z\in x\hplus y$, and $\bar x\,\bar y\coloneqq\overline{xy}$, are well defined.

\begin{lemma}\label{lem:R}
$R$ is a commutative ring. It is local, its maximal ideal is $\bar\varpi R$, and $\bar\varpi^{\,n-1}\ne0=\bar\varpi^{\,n}$. In other words, $R$ is a local Artinian principal ideal ring of length $n$ (a \emph{truncated DVR} in the sense of \cite[Definition~1.8]{Tolliver16}).
\end{lemma}

\begin{proof}
Associativity and commutativity of $+$ follow from (CH1) and (CH2) in $H$, $\bar0$ is neutral because $x\hplus0=\{x\}$, and $\overline{-x}$ is the opposite of $\bar x$ because $0\in x\hminus x$. The multiplication is associative and commutative with unit $\bar1$, and distributive because $xw\in xy\hplus xz$ for $w\in y\hplus z$ by (HF3). Every non-zero $x\in\mathcal O_H$ is $\varpi^j u$ with $j=\val(x)$ and $\val(u)=0$, so $\bar x=\bar\varpi^{\,j}\bar u$, and $\bar x=0$ when $j\ge n$ by \cref{lem:sim}(1). Elements of value $0$ are invertible in $H$ with inverse in $\mathcal O_H$, so $\bar u$ is a unit of $R$, while $\bar\varpi$ is nilpotent: $\bar\varpi^{\,n}=\overline{\varpi^n}=0$, and $\bar\varpi^{\,n-1}\ne0$ because $\val(\varpi^{n-1})<n$. Hence every non-zero element of $R$ is $\bar\varpi^{\,j}\bar u$ with $0\le j<n$ and $\bar u\in R^\times$, so every ideal $I$ of $R$ is $\bar\varpi^{\,j}R$ for the least $j\in\{0,\dots,n\}$ with $\bar\varpi^{\,j}\in I$, and $R$ is local with maximal ideal $\bar\varpi R$.
\end{proof}

\begin{proof}[Proof of \cref{thm:Z}]
By \cref{lem:R}, $R$ is a special principal ideal ring: a principal ideal ring with a unique prime ideal, which is nilpotent. That such a ring is $\mathcal O/\m_{\mathcal O}^n$ for a discrete valuation ring $\mathcal O$ is due to Deligne \cite[p.~126]{Deligne84}, as quoted in the proof of \cite[Theorem~7.5]{Tolliver16}; we use instead a theorem of Hungerford. By \cite[Corollary~11]{Hungerford68}, $R\cong D/I$ for a principal ideal domain $D$ and an ideal $I$. As $R$ is local, $I$ lies in a unique maximal ideal $P$ of $D$, elements of $D\setminus P$ are units modulo $I$, and so $R\cong D_P/ID_P$. Since $n\ge2$, $R$ is not a field, so $P\ne0$ and $\mathcal O\coloneqq D_P$ is a discrete valuation ring with $R\cong\mathcal O/\m_{\mathcal O}^n$, because $R$ has length $n$. This gives a ring isomorphism $\theta\colon\mathcal O/\m_{\mathcal O}^n\to R$. Choose $\pi\in\mathcal O$ with $\theta(\pi\bmod\m^n_{\mathcal O})=\bar\varpi$; since $\bar\varpi\notin(\bar\varpi R)^2$, $\pi$ is a uniformiser. Let $K=\operatorname{Frac}\mathcal O$ with its discrete valuation $v$ and $T=1+\m_\rho(v)=1+\m^n_{\mathcal O}$.

\emph{The map.} Every $a\in K^\times$ is $\pi^k\tilde u$ with $k=v(a)$ and $\tilde u\in\mathcal O^\times$. By \cref{lem:R} and \cref{lem:sim}(2), the unit $\theta(\tilde u\bmod\m^n_{\mathcal O})$ of $R$ is the class of exactly one $u\in H$ of value $0$. Put $\Phi(\cls a_T)\coloneqq \varpi^ku$ and $\Phi(0)=0$. This is well defined, because $\tilde u$ and $\tilde us$ with $s\in T$ have the same residue modulo $\m^n_{\mathcal O}$. It is injective: $\varpi^ku=\varpi^{k'}u'$ forces $k=k'$ and $u=u'$, hence $\tilde u\equiv\tilde u'\bmod\m^n_{\mathcal O}$ and $\tilde u/\tilde u'\in T$. It is surjective, since $x\in H^\times$ equals $\varpi^{\val(x)}(\varpi^{-\val(x)}x)$ and every unit of $R$ has the form $\theta(\tilde u\bmod\m_{\mathcal O}^n)$. It is multiplicative because $\theta$ is, by \cref{lem:sim}(2), and it preserves valuations.

\emph{Compatibility with $R$.} In $\Hr(K)$, with $\cls\pi_T$ in the role of $\varpi$, the construction above gives the ring $\mathcal O/\m^n_{\mathcal O}$. Indeed, for $a,b\in\mathcal O$ every element of $\cls a_T\hminus\cls b_T$ is the class of $a-b-b(s-1)$ with $s\in T$, so $\cls a_T\sim\cls b_T$ if and only if $a\equiv b\bmod\m^n_{\mathcal O}$; and $\cls{a+b}_T\in\cls a_T\hplus\cls b_T$, $\cls{ab}_T=\cls a_T\cls b_T$ show that the induced operations are those of $\mathcal O/\m^n_{\mathcal O}$. For $a=\pi^j\tilde u\in\mathcal O\setminus\{0\}$ we have $\overline{\Phi(\cls a_T)}=\bar\varpi^{\,j}\bar u=\theta(\pi^j\tilde u\bmod\m^n_{\mathcal O})$, and $\overline{\Phi(0)}=\bar 0=\theta(0)$. Since $\Phi$ maps $\mathcal O_{\Hr(K)}$ bijectively onto $\mathcal O_H$ and $\theta$ is bijective, $\Phi$ carries the $\sim$-class of $\cls a_T$ onto the $\sim$-class $\theta(a\bmod\m^n_{\mathcal O})$.

\emph{Sums.} Let $x,y\in\Hr(K)^\times$, and multiply by a power of $\cls\pi_T$ (note $\Phi(\cls\pi_T)=\varpi$) so that $\min\{v(x),v(y)\}=0$; since $\Phi$ is multiplicative and both sides satisfy (HF3), this is harmless. By \cref{lem:sim}(3), applied in $\Hr(K)$ and in $H$, $x\hplus y$ is the class of $\bar x+\bar y$ in $\mathcal O/\m^n_{\mathcal O}$ and $\Phi(x)\hplus\Phi(y)$ is the class of $\overline{\Phi(x)}+\overline{\Phi(y)}=\theta(\bar x+\bar y)$ in $R$. By the previous paragraph, $\Phi(x\hplus y)=\Phi(x)\hplus\Phi(y)$. Sums with $0$ are trivially preserved. Hence $\Phi\colon\Hr(K)\to H$ is an isomorphism of valued hyperfields.
\end{proof}

\begin{remark}\label{rem:discrete-uses}
Discreteness entered in three places: $\val(u)>\rho$ means $\val(u)\ge n$, so that $\m_\rho=\varpi^n\mathcal O_H$; the element $\varpi$ splits $H^\times\to\Z$ and gives the normal form $x=\varpi^ku$; and every truncated DVR is a quotient of a discrete valuation ring. For a general value group and a norm such as $\{\gamma:0\le\gamma\le\delta\}$, the ring $R$ is still defined (it is $\mathcal O_H$ modulo the open ball of radius $\delta$), but one would need to realise $R$ as a quotient of a valuation ring with value group $\Gamma$ and to replace the normal form $\varpi^ku$. We do not know whether either step can fail.
\end{remark}

\begin{remark}[A heuristic for Question~1$'$; not a proof]\label{rem:heuristic}
In Question~1$'$, the norm $\rho$ lies in $\Delta_{\ge0}$ for the convex subgroup $\Delta=\{0\}\times\Z$. Inside $\Delta\cong\Z$, it looks like the norm $\{0,1\}$ of \cref{thm:Z}. One may hope to attach to $H$ a valued hyperfield with value group $\Delta$, apply \cref{thm:Z} to it, and then lift the resulting field along $\Gamma/\Delta\cong\Z$ in the spirit of \cite[Theorem~4.1]{LT25}. We have not carried this out, and we do not know which hyperfield with value group $\Delta$ to attach to $H$.
\end{remark}

\begin{remark}[Tolliver's route]\label{rem:tolliver}
\Cref{thm:Z} can also be obtained from \cite{Tolliver16}. A valued hyperfield of norm $\{0,\dots,n-1\}$ with value group $\Z$ becomes a discretely valued hyperfield in Krasner's sense \cite[Definition~1.4]{Tolliver16} once one puts $\lvert x\rvert=c^{\val(x)}$ for a fixed $0<c<1$. Tolliver attaches to it a triple $\mathrm{Tr}(H)$ in Deligne's sense, whose ring is our $R$. By the proof of \cite[Theorem~7.5]{Tolliver16}, $\mathrm{Tr}(H)\cong\mathrm{Tr}(\Hr(K))$ for $K$ the fraction field of a discrete valuation ring lifting $R$. By \cite[Theorem~1.12]{Tolliver16}, this isomorphism and its inverse lift to morphisms between $H$ and $\Hr(K)$, possibly after rescaling absolute values. An isomorphism of triples has ramification index $r=1$, since the indices multiply; so if the absolute value on $\Hr(K)$ is normalised by the same constant $c$ as on $H$, no rescaling is needed, both lifts are morphisms between the same two objects, and faithfulness \cite[Corollary~6.2]{Tolliver16} makes them mutually inverse. \Cref{thm:A} then upgrades the resulting isomorphism of hyperfields to one of valued hyperfields. We have not written out the dictionary between Krasner's axioms with absolute values and \cref{def:valuation} in detail, which is why we give the direct proof above.
\end{remark}

\section{Massouros hyperfields}\label{sec:massouros}

Massouros~\cite[\S3]{Massouros85} gave the first hyperfields that are not factor hyperfields, answering a question of Krasner~\cite{Krasner83}. We recall his construction and include a proof of his theorem, whose hypothesis is that $G$ has exponent~$2$ and at least four elements; like his, the proof works over any field. For finite $G$ the result is also a special case of a remark of Baker and Jin \cite[Remark~4.4]{BakerJin21}: a finite hyperfield $H$ with $H^\times$ not cyclic and $1\hminus1\ne H$ is not a factor hyperfield.\footnote{Baker and Jin state the remark without finiteness, but its proof uses the theorem of Bergelson--Shapiro and Turnwald \cite[Theorem~4.1]{BakerJin21}, which concerns subgroups of finite index in the multiplicative group of an infinite field. For infinite $H$ the statement fails, for instance for $\mathcal H_{\{0\}}(\mathbb Q_3)$.} Some lower bound on $|G|$ is necessary (\cref{rem:three}).

\begin{definition}[{\cite[\S3]{Massouros85}}]\label{def:massouros}
Let $G$ be an abelian group and $\Mas(G)=G\cup\{0\}$, with the multiplication of $G$ extended by $x\cdot0=0\cdot x=0$. Put $x\hplus 0=0\hplus x=\{x\}$ for all $x\in\Mas(G)$, and for $x,y\in G$
\[
  x\hplus y=\begin{cases}\{x,y\} & \text{if } x\ne y,\\ \Mas(G)\setminus\{x\} & \text{if } x=y.\end{cases}
\]
\end{definition}

In words: two different elements add to the pair they form, and an element added to itself gives everything except itself. In particular $-x=x$ for every $x$. In the running example, $1\hplus1=\{0,a,b,ab\}$ and $a\hplus b=\{a,b\}$ in $\Mas_4$.

\begin{lemma}[{\cite[Proposition~2]{Massouros85}}]\label{lem:massouros-hf}
For every abelian group $G$, $\Mas(G)$ is a hyperfield.
\end{lemma}

\begin{proof}
Commutativity is built in. For (CH3), $0\in x\hplus y$ with $x\ne0$ happens only for $y=x$. For (CH4), let $z\in x\hplus y$; we must show $x\in z\hplus y$. If $x$ or $y$ is $0$ this is immediate. If $x\ne y$, then $z\in\{x,y\}$, and $x\in x\hplus y$ or $x\in y\hplus y=\Mas(G)\setminus\{y\}$. If $x=y$, then $z\ne x$ and $x\in z\hplus x$, which is $\{x\}$ or $\{z,x\}$.

For associativity we may assume $x,y,z\ne0$. If they are pairwise distinct, both $(x\hplus y)\hplus z$ and $x\hplus(y\hplus z)$ equal $\{x,y,z\}$. If exactly two coincide, both sides equal $\Mas(G)$: for instance $(x\hplus x)\hplus z\supseteq (z\hplus z)\cup(0\hplus z)=\Mas(G)$, and $x\hplus(x\hplus z)=(x\hplus x)\cup(x\hplus z)=\Mas(G)$; the other placements are the same computation up to commutativity. If $x=y=z$, the two sides agree by commutativity. Finally, multiplication by $g\in G$ is a bijection of $\Mas(G)$ fixing $0$, and it carries $\{x,y\}$ to $\{gx,gy\}$ and $\Mas(G)\setminus\{x\}$ to $\Mas(G)\setminus\{gx\}$; this is (HF3).
\end{proof}

\begin{proposition}\label{prop:massouros-not-factor}
Let $G$ be an abelian group in which $g^2=1$ for every $g$, and suppose $|G|\ge4$. Then $\Mas(G)$ is not a factor hyperfield.
\end{proposition}

The idea: in a factor hyperfield $k_S$ isomorphic to $\Mas(G)$, the subgroup $S$ would contain $-1$ and all squares, and no element $1+s$ with $s\in S$ could lie in $S$. In a field with many squares one can always write a square as a sum of two elements of $S$.

\begin{proof}
Suppose $\psi\colon k_S\to\Mas(G)$ is an isomorphism, with $k$ a field and $S\le k^\times$. We collect three properties of $S$.
\begin{enumerate}
\item[(a)] $-1\in S$. Indeed $-x=x$ in $\Mas(G)$, so $\cls{-1}_S=-\cls1_S=\cls1_S$.
\item[(b)] $x^2\in S$ for all $x\in k^\times$, because every element of $G$ squares to $1$.
\item[(c)] There is no $s\in S$ with $1+s\in S$. Otherwise $\cls1_S=\cls{1+s}_S\in\cls1_S\hplus\cls1_S$, whereas $1\notin1\hplus1$ in $\Mas(G)$.
\end{enumerate}
Moreover $|k^\times/S|=|G|\ge4$.

We first find $y\in S$ with $y\ne\pm1$. If every square in $k^\times$ were $\pm1$, then every $x\in k^\times$ would satisfy $x^4=1$, so $|k^\times|\le4$. Together with $|k^\times/S|\ge4$ this forces $S=\{1\}$ and $|k^\times|=4$, so $k=\F_5$; but then (a) gives $-1=1$, which is false in $\F_5$. Hence some square $y=x^2$ differs from $\pm1$, and $y\in S$ by (b).

If the characteristic of $k$ is not $2$, put $a=\bigl(\tfrac{y+1}{2}\bigr)^2$ and $b=-\bigl(\tfrac{y-1}{2}\bigr)^2$. Both are non-zero because $y\ne\pm1$, both lie in $S$ by (a) and (b), and $a+b=y\in S$. Hence $s\coloneqq b/a\in S$ and $1+s=y/a\in S$, contradicting (c).

If the characteristic is $2$, then $1+y=(1+x)^2$, which is non-zero since $y\ne1$. So $1+y\in S$ by (b), again contradicting (c).
\end{proof}

\begin{remark}\label{rem:three}
For $|G|=2$ the conclusion fails: $\Mas(\Z/2)$ is isomorphic to $\F_5/\{\pm1\}$, which is case~(2) of the classification of hyperfields of order~$3$ in the proof of \cite[Theorem~1.3(2)]{BakerJin21}. Indeed, in $\F_5/\{\pm1\}=\{0,\cls1,\cls2\}$ we have $\cls1\hplus\cls1=\{\cls2,0\}$, $\cls2\hplus\cls2=\{\cls1,0\}$ and $\cls1\hplus\cls2=\{\cls1,\cls2\}$, which is the table of $\Mas(\Z/2)$. This is checked by machine in \cref{sec:checks}.
\end{remark}

As observed in \cite[\S7(2)]{Linzi25}, a Massouros hyperfield admits only the trivial valuation, even in the weak sense (V0)--(V2): for $x\ne1$ we have $x\in1\hplus1$, so $\val(x)\ge0$, and applying this to $x^{-1}$ gives $\val(x)=0$. To obtain a non-trivial valuation we spread $\Mas(G)$ along a value group.

\section{Layered hyperfields}\label{sec:layer}

Think of an element of a power series field $k((t^\Gamma))$ through its leading term $c\,t^\gamma$. Two leading terms of different degree add to the one of lower degree. Two of the same degree add to $(c+c')\,t^\gamma$, unless the coefficients cancel, in which case anything of higher degree can result. The layering below copies this rule with the field $k$ replaced by an arbitrary hyperfield $\Mas$. This is the $\Gamma$-layering of Bowler and Su \cite[Definition~4.1]{BowlerSu21} along the split exact sequence of groups $1\to\Mas^\times\to\Mas^\times\times\Gamma\to\Gamma\to1$ (the last two written multiplicatively, as in \cite{BowlerSu21}), written with the opposite order on $\Gamma$ (Bowler and Su let larger elements absorb smaller ones); see also the proof of \cite[Proposition~2.5]{LT25}. That it is a hyperfield is \cite[Lemma~4.2]{BowlerSu21}, whose additive part rests on \cite[Lemma~3.1]{BowlerSu21}, valid for every hyperfield $\Mas$. To keep the note self-contained we include a short proof in our notation, organised around the absorption formula of \cref{lem:union}.

Fix a hyperfield $\Mas$ and an ordered abelian group $\Gamma$. As a set, $\Lay=\Mas\rtimes\Gamma$ consists of $0$ and the pairs $(c,\gamma)$ with $c\in\Mas^\times$ and $\gamma\in\Gamma$; we write $c\,t^\gamma$ for $(c,\gamma)$ and $t^\gamma$ for $1\,t^\gamma$. The unit of $\Lay$ is $1=t^0$. Multiplication is $(c\,t^\gamma)(c'\,t^{\gamma'})=cc'\,t^{\gamma+\gamma'}$, and $0$ is absorbing. The map $v(c\,t^\gamma)=\gamma$, $v(0)=\infty$, will be the valuation.

\begin{definition}\label{def:At}
For $\gamma\in\Gamma$ let $B_\gamma\coloneqq\{x\in\Lay:v(x)>\gamma\}$; it contains $0$. For a non-empty subset $A\subseteq\Mas$ put
\[
  \lt{A}{\gamma}\coloneqq\{c\,t^\gamma: c\in A\setminus\{0\}\}\cup
  \begin{cases} B_\gamma & \text{if } 0\in A,\\ \varnothing & \text{if } 0\notin A.\end{cases}
\]
\end{definition}

So $\lt{A}{\gamma}$ is ``$A$ placed in degree $\gamma$'', where a cancellation in $A$ is replaced by everything of higher degree. Unions behave well: $\bigl(\bigcup_iA_i\bigr)t^\gamma=\bigcup_i\lt{A_i}{\gamma}$.

\begin{definition}\label{def:layer}
The hyperaddition of $\Lay$ is $x\hplus0=0\hplus x=\{x\}$ and, for $x=c\,t^\gamma$ and $y=c'\,t^{\gamma'}$,
\[
  x\hplus y\coloneqq\begin{cases}
  \{x\} & \text{if }\gamma<\gamma',\\
  \{y\} & \text{if }\gamma'<\gamma,\\
  \lt{(c\hplus c')}{\gamma} & \text{if }\gamma=\gamma'.
  \end{cases}
\]
\end{definition}

A uniform description is convenient. For non-zero $x_1,\dots,x_n\in\Lay$ with $x_i=c_i\,t^{\gamma_i}$, let $\mu=\min_i\gamma_i$ and let $\sigma(x_1,\dots,x_n)\subseteq\Mas$ be the hypersum, computed in $\Mas$, of the $c_i$ with $\gamma_i=\mu$ (the singleton $\{c_i\}$ if only one index attains $\mu$). Then, for non-zero $x,y$,
\begin{equation}\label{eq:two-sum}
  x\hplus y=\lt{\sigma(x,y)}{\mu},\qquad \mu=\min\{v(x),v(y)\}.
\end{equation}

\begin{example}\label{ex:L4}
In $\Lay_4=\Mas_4\rtimes\Z$ we have $t^0\hplus a\,t^1=\{t^0\}$ and $a\,t^0\hplus b\,t^0=\{a\,t^0,b\,t^0\}$, while
\[
  t^0\hplus t^0=\{a\,t^0, b\,t^0, ab\,t^0\}\cup B_0,
\]
since $1\hplus1=\{0,a,b,ab\}$ in $\Mas_4$ contains $0$. So the sum $1\hplus1=1\hminus1$ in $\Lay_4$ contains elements of value $0$, whereas \cref{lem:one-minus-one} requires every element of $1\hminus1$ in a valued hyperfield to have value $>0$ (see \cref{lem:not-krasner}).
\end{example}

The key computation is how a set of the form $\lt{A}{\alpha}$ absorbs one more summand.

\begin{lemma}\label{lem:union}
Let $\varnothing\ne A\subseteq\Mas$, $\alpha\in\Gamma$ and $z=c\,t^\gamma\in\Lay^\times$. Then
\[
  \lt{A}{\alpha}\hplus z=\begin{cases}
  \lt{A}{\alpha} & \text{if }\alpha<\gamma,\\
  \lt{\{c\}}{\gamma} & \text{if }\gamma<\alpha,\\
  \lt{(A\hplus c)}{\gamma} & \text{if }\alpha=\gamma.
  \end{cases}
\]
\end{lemma}

\begin{proof}
\emph{Case $\gamma<\alpha$.} Every element of $\lt{A}{\alpha}$ has value $>\gamma$, so its sum with $z$ is $\{z\}$.

\emph{Case $\alpha<\gamma$.} An element $d\,t^\alpha$ with $d\in A\setminus\{0\}$ contributes $\{d\,t^\alpha\}$. If $0\in A$ we must also add $z$ to every $w\in B_\alpha$, and we claim that these sums fill exactly $B_\alpha$. Each sum lies in $B_\alpha$ by \cref{def:layer}. Conversely, $w=0$ gives $z$; each $w$ with $\alpha<v(w)<\gamma$ gives $w$; each $w$ with $v(w)>\gamma$ gives $z$ again; the elements $w=d\,t^\gamma$ give $\bigcup_{d\in\Mas^\times}\lt{(d\hplus c)}{\gamma}$, which contains $c'\,t^\gamma$ for every $c'\ne 0$ (as $\bigcup_{d\in\Mas}d\hplus c=\Mas$ and $d=0$ gives $c$ itself) and contains $B_\gamma$ (take $d=-c$). Hence the union is $B_\alpha$.

\emph{Case $\alpha=\gamma$.} An element $d\,t^\gamma$ with $d\in A\setminus\{0\}$ contributes $\lt{(d\hplus c)}{\gamma}$. If $0\in A$, the elements of $B_\gamma$ contribute only $z=\lt{\{c\}}{\gamma}=\lt{(0\hplus c)}{\gamma}$, since everything in $B_\gamma\setminus\{0\}$ has value $>\gamma$. Taking the union over $d\in A$ gives $\lt{(A\hplus c)}{\gamma}$.
\end{proof}

\begin{proposition}\label{prop:layer-hf}
For every hyperfield $\Mas$ and every ordered abelian group $\Gamma$, the layering $\Lay=\Mas\rtimes\Gamma$ is a hyperfield, and $v$ satisfies \textup{(V0)--(V2)}.
\end{proposition}

\begin{proof}
\emph{Associativity.} If one of $x,y,z$ is $0$, both sides of (CH1) reduce to a sum of the other two. Otherwise, \eqref{eq:two-sum} and \cref{lem:union} give
\[
 (x\hplus y)\hplus z=\lt{\sigma(x,y,z)}{\mu},\qquad \mu=\min\{v(x),v(y),v(z)\}.
\]
To see this, apply the lemma with $A=\sigma(x,y)$ and $\alpha=\min\{v(x),v(y)\}$. If $\alpha<v(z)$, then $z$ is not among the summands of least value. If $v(z)<\alpha$, it is the only one. If $\alpha=v(z)$, we use $\sigma(x,y)\hplus c=\sigma(x,y,z)$, which is associativity in $\Mas$. The right-hand side is symmetric in $x,y,z$, so by (CH2) it also equals $(y\hplus z)\hplus x=x\hplus(y\hplus z)$.

\emph{Commutativity} is clear. \emph{Inverses:} for non-zero $x=c\,t^\gamma$ and $y$, we have $0\in x\hplus y$ if and only if $v(y)=\gamma$ and $0\in c\hplus c'$, where $y=c'\,t^\gamma$; that is, if and only if $y=(-c)\,t^\gamma$. So $-x=(-c)\,t^\gamma$ is unique.

\emph{Reversibility.} Let $z\in x\hplus y$; we show $x\in z\hminus y$. The cases with a zero among $x,y,z$ are immediate (if $z=0$ then $y=-x$ and $x\in 0\hplus x$). Let $x=c\,t^\alpha$ and $y=c'\,t^\beta$. If $\alpha<\beta$, then $z=x$ and $x\hminus y=\{x\}$. If $\beta<\alpha$, then $z=y$ and $x\in B_\beta\subseteq y\hminus y$. If $\alpha=\beta$, either $z=e\,t^\alpha$ with $e\in c\hplus c'$, in which case $c\in e\hminus c'$ by (CH4) in $\Mas$ and so $x\in z\hminus y$; or $c'=-c$ and $v(z)>\alpha$, in which case $z\hminus y=z\hplus x=\{x\}$.

\emph{Multiplication and distributivity.} $\Lay^\times=\Mas^\times\times\Gamma$ is an abelian group and $0$ is absorbing. Multiplication by $d\,t^\delta$ maps $\lt{A}{\gamma}$ onto $\lt{(dA)}{\gamma+\delta}$, because it maps $B_\gamma$ onto $B_{\gamma+\delta}$ and $0\in dA$ if and only if $0\in A$. Since $d\,\sigma(x,y)=\sigma(dx,dy)$ by (HF3) in $\Mas$, this gives (HF3) in $\Lay$.

\emph{Valuation.} (V0) and (V1) are clear and $v$ is surjective. Every element of $\lt{A}{\mu}$ has value at least $\mu$, which is (V2).
\end{proof}

\section{Proof of Theorem~\ref{thm:transfer}}\label{sec:B}

Throughout this section $\Mas$ is a hyperfield, $\Gamma$ an ordered abelian group and $\Lay=\Mas\rtimes\Gamma$. By \cref{prop:layer-hf}, $\Lay$ is a hyperfield and $v$ satisfies (V0)--(V2); it is surjective, and non-trivial when $\Gamma\ne0$. We prove the two directions of the equivalence separately (\cref{prop:residue,prop:lift}) and then specialise to Massouros hyperfields.

An isomorphism of hyperfields $f\colon\Mas\to\Mas'$ induces an isomorphism $\Mas\rtimes\Gamma\to\Mas'\rtimes\Gamma$, $c\,t^\gamma\mapsto f(c)\,t^\gamma$, because $f$ maps each $\lt{A}{\gamma}$ onto $\lt{f(A)}{\gamma}$ and $0\in f(A)$ if and only if $0\in A$. So being a factor hyperfield, for $\Mas$ and for $\Lay$, depends only on the isomorphism class of $\Mas$.

\subsection{From the layering to the residue}
A factor hyperfield with the layered structure of $\Lay$ has $\Mas$ as its residue. In the language of \cite{Linzi25}, this is \cite[Lemma~3.6]{Linzi25}, which lifts the valuation, combined with \cite[Proposition~3.14]{Linzi25}, which identifies the residue hyperfield of $K_T$ with the factor hyperfield $(Kw)_{\bar T}$ of the residue field; the residue hyperfield of $\Lay$ is visibly $\Mas$. We give a direct proof.

\begin{proposition}\label{prop:residue}
If $\Lay$ is isomorphic to a factor hyperfield $K_T$, then $\Mas$ is isomorphic to a factor hyperfield of a field, namely of the residue field of a valuation on $K$.
\end{proposition}

\begin{proof}
Let $\varphi\colon K_T\to\Lay$ be an isomorphism and put $w(x)\coloneqq v(\varphi(\cls x_T))$. Step~1 of the proof of \cref{thm:A} uses only (V0)--(V2), so $w$ is a valuation on $K$. Let $\mathcal O_w$, $\m_w$ and $k$ be its valuation ring, maximal ideal and residue field, and let $\res\colon\mathcal O_w\to k$ be the residue map.

\emph{$T$ lies between $1+\m_w$ and $\mathcal O_w^\times$.} Every $s\in T$ has $w(s)=v(1)=0$. If $w(m)>0$, then $1\hplus\varphi(\cls m_T)=\{1\}$ by \cref{def:layer}, so $\cls{1+m}_T=\cls1_T$, since $\cls{1+m}_T\in\cls1_T\hplus\cls m_T$; hence $1+m\in T$.

Let $\bar T=\res(T)\le k^\times$ and let $k_{\bar T}$ be the factor hyperfield. Elements of value $0$ in $\Lay$ are the $c\,t^0$ with $c\in\Mas^\times$. For such $c$ choose $u_c\in\mathcal O_w^\times$ with $\varphi(\cls{u_c}_T)=c\,t^0$, and define
\[
  \psi\colon\Mas\to k_{\bar T},\qquad \psi(0)=0,\quad \psi(c)=\cls{\res(u_c)}_{\bar T}.
\]

\emph{$\psi$ is a well-defined isomorphism of groups $\Mas^\times\to k^\times/\bar T$.} Another choice of $u_c$ differs by a factor in $T$, so $\psi$ is well defined; it is multiplicative because $\varphi$ is. It is surjective because $\res(\mathcal O_w^\times)=k^\times$. If $\res(u_c)=\res(u_{c'}s)$ with $s\in T$, then $u_c/(u_{c'}s)\in1+\m_w\subseteq T$, so $\cls{u_c}_T=\cls{u_{c'}}_T$ and $c=c'$.

\emph{$\psi$ preserves sums.} Let $c,c'\in\Mas^\times$ and write $u=u_c$, $u'=u_{c'}$. By \cref{def:layer} and \eqref{eq:factor},
\[
  \lt{(c\hplus c')}{0}=c\,t^0\hplus c'\,t^0=\{\varphi(\cls{u+u's}_T):s\in T\}.
\]
Recall that $\psi(c)\hplus\psi(c')=\{\cls{\res(u)+\res(u')\tau}_{\bar T}:\tau\in\bar T\}$. Let $e\in\Mas^\times$. If $e\in c\hplus c'$, then $e\,t^0=\varphi(\cls{u+u's}_T)$ for some $s\in T$; thus $w(u+u's)=0$, we may take $u_e=u+u's$, and $\psi(e)=\cls{\res(u)+\res(u')\res(s)}_{\bar T}\in\psi(c)\hplus\psi(c')$. Conversely, suppose $\psi(e)=\cls{\res(u)+\res(u')\res(s)}_{\bar T}\ne0$ for some $s\in T$. Then $w(u+u's)=0$, so $\varphi(\cls{u+u's}_T)=e'\,t^0$ for some $e'\in c\hplus c'$, and $\psi(e')=\psi(e)$ by the first part. Injectivity gives $e=e'\in c\hplus c'$. For $e=0$: $0\in c\hplus c'$ if and only if $\lt{(c\hplus c')}{0}$ contains $0$ or an element of positive value, that is, if and only if $w(u+u's)>0$ for some $s\in T$. This in turn means $\res(u)+\res(u')\res(s)=0$, that is, $0\in\psi(c)\hplus\psi(c')$. Sums involving $0$ are trivially preserved. So $\psi$ is an isomorphism of hyperfields.
\end{proof}

\subsection{From the residue to the layering}
For the converse we need a valued field with residue field $k$, value group $\Gamma$ and a multiplicative leading coefficient. For $S=k^\times$, for $S$ the positive cone of an ordered field and for $S=\{1\}$, Bowler and Su \cite[proof of Theorem~7.5, Cases~1--3]{BowlerSu21} use Hahn-type fields and the same subgroup $T$ as below; we extend their construction to arbitrary $S$. The Hahn field $k(\!(t^\Gamma)\!)$ would do; the following smaller field is enough. Let $k[t^\Gamma]$ be the group algebra of $\Gamma$ over $k$: its elements are finite sums $f=\sum_\gamma c_\gamma t^\gamma$ with $c_\gamma\in k$, multiplied by $t^\gamma t^\delta=t^{\gamma+\delta}$. For $f\ne0$ let $v(f)$ be the least $\gamma$ with $c_\gamma\ne0$ and $\operatorname{lc}(f)\coloneqq c_{v(f)}$.

\begin{lemma}\label{lem:group-field}
$k[t^\Gamma]$ is an integral domain. Let $k(t^\Gamma)$ be its field of fractions. The maps $v$ and $\operatorname{lc}$ extend uniquely to group homomorphisms $v\colon k(t^\Gamma)^\times\to\Gamma$ and $\operatorname{lc}\colon k(t^\Gamma)^\times\to k^\times$, with $v(c\,t^\gamma)=\gamma$ and $\operatorname{lc}(c\,t^\gamma)=c$. For $x,y\in k(t^\Gamma)^\times$:
\begin{enumerate}
\item if $v(x)<v(y)$, then $v(x+y)=v(x)$ and $\operatorname{lc}(x+y)=\operatorname{lc}(x)$;
\item if $v(x)=v(y)$ and $\operatorname{lc}(x)+\operatorname{lc}(y)\ne0$, then $v(x+y)=v(x)$ and $\operatorname{lc}(x+y)=\operatorname{lc}(x)+\operatorname{lc}(y)$;
\item if $v(x)=v(y)$ and $\operatorname{lc}(x)+\operatorname{lc}(y)=0$, then $x+y=0$ or $v(x+y)>v(x)$.
\end{enumerate}
In particular $v$, with $v(0)=\infty$, is a valuation on $k(t^\Gamma)$ with value group $\Gamma$.
\end{lemma}

\begin{proof}
For non-zero $f,g\in k[t^\Gamma]$, every term of $fg$ other than $\operatorname{lc}(f)\operatorname{lc}(g)\,t^{v(f)+v(g)}$ has degree $\gamma+\delta$ with $\gamma\ge v(f)$, $\delta\ge v(g)$ and one inequality strict, so of degree $>v(f)+v(g)$, because $\Gamma$ is ordered. Hence $fg\ne0$, $v(fg)=v(f)+v(g)$ and $\operatorname{lc}(fg)=\operatorname{lc}(f)\operatorname{lc}(g)$. So $v(f/g)\coloneqq v(f)-v(g)$ and $\operatorname{lc}(f/g)\coloneqq\operatorname{lc}(f)/\operatorname{lc}(g)$ are well defined homomorphisms, and they are the only extensions. For (1)--(3), write $x=f/g$ and $y=h/g$ with a common denominator. Then $x+y=(f+h)/g$, and the lowest term of $f+h$ is that of $f$ if $v(f)<v(h)$, and is $(\operatorname{lc}(f)+\operatorname{lc}(h))\,t^{v(f)}$ if $v(f)=v(h)$ and this coefficient is non-zero; if it is zero, $f+h=0$ or $v(f+h)>v(f)$. Dividing by $g$ gives (1)--(3), and these imply $v(x+y)\ge\min\{v(x),v(y)\}$.
\end{proof}

\begin{proposition}\label{prop:lift}
Let $k$ be a field, $S\le k^\times$ a subgroup, $K=k(t^\Gamma)$ and
\[
  T\coloneqq\{x\in K^\times: v(x)=0,\ \operatorname{lc}(x)\in S\}.
\]
Then $\Phi\colon K_T\to k_S\rtimes\Gamma$, $\Phi(\cls x_T)=\cls{\operatorname{lc}(x)}_S\,t^{v(x)}$, $\Phi(0)=0$, is an isomorphism of hyperfields. Consequently, if $\Mas$ is a factor hyperfield, so is $\Mas\rtimes\Gamma$.
\end{proposition}

\begin{proof}
$T$ is a subgroup of $K^\times$, because $v$ and $\operatorname{lc}$ are homomorphisms. For $x,y\in K^\times$ we have $\Phi(\cls x_T)=\Phi(\cls y_T)$ if and only if $v(x/y)=0$ and $\operatorname{lc}(x/y)\in S$, that is, $x/y\in T$; so $\Phi$ is well defined and injective. It is surjective because $\Phi(\cls{c\,t^\gamma}_T)=\cls c_S\,t^\gamma$, and it is multiplicative. It remains to show $\Phi(\cls x_T\hplus\cls y_T)=\Phi(\cls x_T)\hplus\Phi(\cls y_T)$ for $x,y\in K^\times$; sums with $0$ are clear. Put $\alpha=v(x)$, $\beta=v(y)$, $a=\operatorname{lc}(x)$, $b=\operatorname{lc}(y)$, and recall that $\cls x_T\hplus\cls y_T=\{\cls{x+ys}_T:s\in T\}$.

If $\alpha<\beta$, then $v(ys)=\beta$ for every $s\in T$, so $x+ys$ has value $\alpha$ and leading coefficient $a$ by \cref{lem:group-field}(1). Both sides equal $\{\cls a_S\,t^\alpha\}$ by \cref{def:layer}. The case $\beta<\alpha$ is symmetric.

Let $\alpha=\beta$. By \cref{def:layer} and \eqref{eq:factor} for $k_S$, the right-hand side is $\lt{(\cls a_S\hplus\cls b_S)}{\alpha}$, where $\cls a_S\hplus\cls b_S=\{\cls{a+b\lambda}_S:\lambda\in S\}$. For $s\in T$ put $\lambda=\operatorname{lc}(s)\in S$; then $v(ys)=\alpha$ and $\operatorname{lc}(ys)=b\lambda$. If $a+b\lambda\ne0$, then $\Phi(\cls{x+ys}_T)=\cls{a+b\lambda}_S\,t^\alpha$ by \cref{lem:group-field}(2). If $a+b\lambda=0$, then $\Phi(\cls{x+ys}_T)\in B_\alpha$ by \cref{lem:group-field}(3). This proves $\subseteq$. Conversely, let $\lambda\in S$. Taking $s=\lambda\in T$ gives $\cls{a+b\lambda}_S\,t^\alpha$ whenever $a+b\lambda\ne0$. Suppose $a+b\lambda=0$; we show that $B_\alpha$ lies in the left-hand side. Let $z\in K$ with $z=0$ or $v(z)>\alpha$. By \cref{lem:group-field}(3), $x+y\lambda$ is $0$ or has value $>\alpha$, so $m\coloneqq(z-x-y\lambda)/y$ is $0$ or has value $>0$. Then $s\coloneqq\lambda+m$ has $v(s)=0$ and $\operatorname{lc}(s)=\lambda$ by \cref{lem:group-field}(1), so $s\in T$, and $x+ys=z$. As $\Phi$ is surjective and preserves values, $\{\Phi(\cls z_T):z=0\text{ or }v(z)>\alpha\}=B_\alpha$. This proves $\supseteq$.

For the last assertion, if $\Mas\cong k_S$, then $\Mas\rtimes\Gamma\cong k_S\rtimes\Gamma\cong K_T$.
\end{proof}

For $S=\{1\}$ we get $T=\{x:v(x-1)>0\}=1+\m_{\{0\}}$, so $k\rtimes\Gamma\cong\mathcal H_{\{0\}}(k(t^\Gamma))$, the first row of \cref{tab:examples}.

\subsection{Massouros hyperfields}

\begin{lemma}\label{lem:not-krasner}
Let $G$ be an abelian group with $|G|\ge2$. Then $\Mas(G)\rtimes\Gamma$ admits no valuation of any norm in the sense of \cref{def:valuation}.
\end{lemma}

The obstruction is the one noted in \cite[\S2, after Theorem~2.2]{LT25} for layerings of the Krasner and sign hyperfields, and behind \cite[Proposition~4.14]{Linzi25}: by \cref{lem:one-minus-one}, $1\hminus1$ may contain no element of value~$0$.

\begin{proof}
Write $\Mas=\Mas(G)$ and $\Lay=\Mas\rtimes\Gamma$, and let $\val'$ be such a valuation, with value group $\Gamma'$ and norm $\rho'$. Since $-1=1$ in $\Mas$, also $-t^0=t^0$ in $\Lay$, so $1\hminus1=1\hplus1=\lt{(1\hplus1)}{0}$, where $1\hplus1=\Mas\setminus\{1\}$ is computed in $\Mas$ (\cref{def:massouros}). Hence $1\hminus1$ contains $g\,t^0$ and $g^{-1}t^0$ for any $g\in G\setminus\{1\}$; see \cref{ex:L4} for $G=V$. By \cref{lem:one-minus-one}, both have value $>\rho'\ni0$ under $\val'$, so $\val'(g\,t^0)>0$ and $\val'(g^{-1}t^0)>0$. But their product is $t^0=1$, so the two values add up to $\val'(1)=0$. This is a contradiction.
\end{proof}

\begin{proof}[Proof of \cref{thm:transfer}]
By \cref{prop:layer-hf}, $\Lay$ is a hyperfield and $v$ is a surjective map satisfying (V0)--(V2). By \cref{prop:residue,prop:lift}, $\Lay$ is a factor hyperfield if and only if $\Mas$ is. If $\Mas=\Mas(G)$ with $G$ of exponent~$2$ and $|G|\ge4$, then $\Mas$ is not a factor hyperfield by \cref{prop:massouros-not-factor}, hence neither is $\Lay$, and \cref{lem:not-krasner} shows that $\Lay$ admits no valuation of any norm.
\end{proof}

By \cref{rem:three}, $\Mas(\Z/2)\cong\F_5/\{\pm1\}$ is a factor hyperfield, so $\Mas(\Z/2)\rtimes\Gamma$ is one too, by \cref{prop:lift}. By \cref{lem:not-krasner} it carries no valuation of any norm. This is another illustration of \cref{rem:audit}: a factor hyperfield with a non-trivial map satisfying (V0)--(V2) need not be any $\Hr(K)$.

\begin{remark}[The rows of \cref{tab:examples}]\label{rem:table}
The Krasner hyperfield $\{0,1\}$, with $1\hplus1=\{0,1\}$, is $k_{k^\times}$ for any field $k$ with at least three elements \cite[Example~2.10]{BakerBowler19}, \cite[Proposition~2.17(i)]{Linzi25}. The hyperfield of signs $\{0,1,-1\}$, with $1\hplus1=\{1\}$ and $1\hplus(-1)=\{0,1,-1\}$, is $\mathbb R_{\mathbb R_{>0}}$ \cite[Example~2.13]{BakerBowler19}, \cite[Proposition~2.17(ii)]{Linzi25}. By \cref{prop:lift}, their layerings are factor hyperfields; this is \cite[proof of Theorem~7.5, Cases~1 and~2]{BowlerSu21}, and for the Krasner layering, which is the generalised tropical hyperfield $\mathcal T(\Gamma)$ of \cite[Example~2.14]{Linzi25}, also \cite[Proposition~2.17(iv)]{Linzi25}. In both layerings $-1=(-1)\,t^0$ and $1\in1\hplus(-1)$ in the coefficient hyperfield, so $1=t^0\in1\hminus1$. By \cref{lem:one-minus-one}, a valuation of some norm would give $\val'(1)>0$, which is absurd. In the Krasner layering, $t^\alpha\hplus t^\beta=\{t^{\min\{\alpha,\beta\}}\}$ for $\alpha\ne\beta$, and $t^\alpha\hplus t^\alpha$ is the set of all elements of value at least $\alpha$. For $\Gamma=\mathbb R$, the bijection $t^\gamma\mapsto-\gamma$, $0\mapsto-\infty$, is therefore an isomorphism onto Viro's tropical hyperfield, in which $a\hplus b=\{\max\{a,b\}\}$ for $a\ne b$ and $a\hplus a=\{c:c\le a\}$ \cite[\S5.3]{Viro10}. Likewise $\pm t^\gamma\mapsto\pm e^{-\gamma}$, $0\mapsto0$, identifies the layering of the hyperfield of signs over $\mathbb R$ with Viro's tropical real hyperfield on $\mathbb R$ \cite[\S7.2]{Viro10}, in which $a\hplus b$ is $\{a\}$ if $|a|>|b|$ or $a=b$, and is $[-|a|,|a|]$ if $a=-b$. For general $\Gamma$ we use the names in \cref{tab:examples} for these layerings. The last row is \cref{thm:transfer}.
\end{remark}

\begin{remark}\label{rem:B-scope}
\Cref{thm:transfer} answers the question of \cite[\S7(2)]{Linzi25} in the negative: for $\Mas=\Mas(G)$ as above and $\Gamma\ne0$, the map $v$ on $\Lay$ is non-trivial and satisfies (V0)--(V2), the axioms of \cite[Definition~3.1]{Linzi25} (labelled (V1)--(V3) there), yet $\Lay$ is not a factor hyperfield. More generally, every hyperfield $\Mas$ that is not a factor hyperfield gives such an example $\Mas\rtimes\Gamma$. This does not bear on the question of Touchard and the author, because $\Lay$ is not a valued hyperfield in the sense of \cref{def:valuation}. It shows that (V3) cannot be dropped from the hypothesis of the question in \cite[Remark~4.12]{Linzi25}, that is, of \cref{cor:equivalent}(2): without the norm axiom, a hyperfield with a non-trivial valuation need not be a factor hyperfield. The norm axiom is also needed in \cref{thm:A} itself; see \cref{rem:audit}. For a Krasner valuation, \cref{lem:one-minus-one} forces the residue structure to be a field \cite[Proposition~4.14]{Linzi25}, so this particular obstruction disappears.
\end{remark}

\section{Computer checks}\label{sec:checks}

The repository accompanying this paper contains short Python scripts that test the finite and combinatorial claims by brute force.
\begin{itemize}
\item \texttt{massouros.py} verifies all hyperfield axioms for $\Mas((\Z/2)^k)$, $k=1,2,3$, and verifies the isomorphism $\Mas(\Z/2)\cong\F_5/\{\pm1\}$ of \cref{rem:three}.
\item \texttt{layering.py} represents subsets $\lt{A}{\gamma}$ symbolically and checks, for $\Mas$ among $\Mas(\Z/2)$, $\Mas_4$, $\F_3$, the Krasner hyperfield and the sign hyperfield, the axioms (CH1)--(CH4) and (HF3), the formula of \cref{lem:union} against an elementwise union, and (V0)--(V2) for $\Mas\rtimes\Z$ on all triples of values in $\{0,1,2\}$. Only the relative order of the three values matters for associativity, so this window covers every case. It also confirms that (V3) fails on $1\hminus1$ for each of the four tested $\Mas$ that is not a field. The script \texttt{mutants.py} confirms that two deliberately broken variants of the construction are rejected by the same checks.
\item \texttt{discrete\_ring.py} checks, in $\Hr(\mathbb Q_p)$ for $\rho=\{0,\dots,n-1\}$ and small $p,n$, that the relation $\sim$ of \cref{lem:sim}, computed from the hyperaddition, is congruence modulo $p^n$, and that sums of elements of minimal value $0$ are full $\sim$-classes. This exercises only the canonical examples.
\item \texttt{factor\_trivial\_valuation.py} checks the smallest instance of \cref{thm:A}: for primes $p\le23$, the factor hyperfield $\F_p/T$ satisfies the norm axiom for the trivial valuation exactly when $T=\{1\}$.
\item \texttt{norm\_axiom\_needed\_for\_thmA.py} checks the first example of \cref{rem:audit}: $\F_5/\{\pm1\}$ with the trivial valuation satisfies (V0)--(V2), (V3) fails, and $T\ne1+\m_\rho$.
\item \texttt{norm\_audit.py} checks the formula of \cref{rem:norm-intrinsic} in $\Hr(\mathbb Q_p)$ for small $p$ and $\rho=\{0,\dots,n-1\}$, and both examples of \cref{rem:audit}.
\item \texttt{factor\_layering.py} checks \cref{prop:lift} for $\Gamma=\Z$ and $k=\F_p$, for six pairs $(p,S)$. Among them are the Krasner hyperfield, $\F_5/\{\pm1\}\cong\Mas(\Z/2)$ and the field $\F_5$. In Laurent series truncated at relative precision $2$ or $3$, it computes the $T$-orbits on units and compares them with the fibres of $u\mapsto\cls{\operatorname{lc}(u)}_S$. For all $x,y$ with values in a window, it computes $\Phi(\cls x_T\hplus\cls y_T)$ and compares it with the layered sum. A deliberately wrong choice of $T$ is rejected. The script also checks the identifications of \cref{rem:table} with Viro's hyperfields on a finite window of $\Z$, and that $1\in1\hminus1$ in the Krasner and sign layerings.
\end{itemize}
These scripts support, but do not replace, the proofs above.

\subsection*{Data availability}
The Python scripts described in this section are contained in the repository accompanying this paper. No other data were generated or analysed.

\subsection*{Acknowledgements}
The author used AI-based tools to assist with drafting, editing and checking this manuscript. All mathematical content was verified by the author, who takes full responsibility for it.

\bibliographystyle{amsplain}
\bibliography{refs}

\end{document}